\documentclass[11pt]{article}
\usepackage{amsfonts}

\usepackage{CJK}
\usepackage{amsfonts,amssymb,amsmath,mathrsfs,multirow,booktabs}
\usepackage{xcolor,soul}
\usepackage{color,latexsym,amsfonts}
 \newcommand{\qed}{\hfill\rule{2mm}{3mm}\vspace{4mm}}

 \newtheorem{theorem}{Theorem}[section]
 \newtheorem{lemma}[theorem]{Lemma}
 \newtheorem{corollary}[theorem]{Corollary}
 \newtheorem{proposition}[theorem]{Proposition}
 \newtheorem{example}[theorem]{Example}
 \newtheorem{Definition}[theorem]{Definition}
 \newtheorem{remark}[theorem]{Remark}
 \newtheorem{condition}[theorem]{Condition}
 \newtheorem{conjecture}[theorem]{Conjecture}

 \def\beqlb{\begin{eqnarray}}\def\eeqlb{\end{eqnarray}}
 \def\beqnn{\begin{eqnarray*}}\def\eeqnn{\end{eqnarray*}}

 \def\<{\langle}\def\>{\rangle}

 \def\ar{&\!\!}

 \def\eqref#1{{\rm(\ref{#1})}}

 \def\proof{\noindent{\it
 Proof.~}}\def\qed{\hfill$\Box$\medskip}

\def\e{{\mbox{\rm e}}}
\def\<{\left<}\def\>{\right>}

\newcommand{\dd}{\mathrm{d}}

\newfam\msbmfam\font\tenmsbm=msbm10\textfont
\msbmfam=\tenmsbm\font\sevenmsbm=msbm7
\scriptfont\msbmfam=\sevenmsbm

\def\<{\left<}\def\>{\right>}
\def\({\left(}\def\){\right)}

\begin{document}

 \centerline{\Large\bf Boundary behaviors for a continuous-state}
\smallskip

 \centerline{\Large\bf nonlinear Neveu's branching process}

   \smallskip\smallskip
 \centerline{Linyu Bai
and Xu Yang\footnote{Corresponding author.}}

    \smallskip\smallskip

 \centerline{School of Mathematics and Information Science,
North Minzu University}

 \centerline{Yinchuan 750021, People's Republic of China}

\smallskip

 \centerline{Emails: bly202205@163.com, xuyang@mail.bnu.edu.cn}

\bigskip

 {\narrower{\narrower

\noindent{\bf Abstract}.
After generalizing the criteria introduced by Chen, we establish
the necessary and sufficient conditions
for the extinction, explosion and coming down from infinity of
a continuous-state nonlinear Neveu's branching  process.

\medskip


\medskip


\par}\par}

\section{Introduction and main result}

\setcounter{equation}{0}

Branching processes in discrete-state are
introduced as probabilistic models for
the stochastic evolution of populations,
and their scaling limits are regarded as
continuous-state branching processes (CSBPs for short),
which are nonnegative-valued
Markov processes with the additive branching processes.
We refer to Kyprianou \cite{Kyprianou7} and Li \cite{Li11,Li19}
and references therein
for reviews and literature on CSBPs.

Recently,
generalized versions of CSBPs have been introduced to study the
interactions between individuals and (or) between individuals and the population.
In Li \cite{Li9}, a class of CSBPs with nonlinear branching mechanism is given
by generalized Lamperti transform.
A more general version of the nonlinear
CSBP is studied in Ma et al. \cite{Ma et al.13} as the solution to
stochastic differential equation (SDE for short)
 \beqlb\label{1.01}
X_t
 \ar=\ar
x+\int_0^t a_0(X_s)\dd s+\int_0^t\int_0^{a_1(X_s)}W(\dd s,\dd u) \cr
 \ar\ar
+\int_0^t\int_U\int_0^{a_2(X_{s-})}z\tilde{N}(\dd s,\dd z,\dd u)
+\int_0^t\int_{(0,\infty)\setminus U}\int_0^{a_3(X_{s-})}zN(\dd s,\dd z,\dd u),
 \eeqlb
where $x>0$, $a_0$ and $a_1,a_2,a_3\ge0$ are Borel function on $[0, \infty)$,
$W(\dd s,\dd u)$ and
$N(\dd s,\dd z,\dd u)$
denote a Gaussian white noise and an
independent Poisson random measure,
respectively. Here
$U$ is a Borel set on $(0,\infty)$ and $\tilde{N}(\dd s,\dd z,\dd u)$
is the compensated measure of $N(\dd s,\dd z,\dd u)$.
The pathwise unique nonnegative solution to SDE \eqref{1.01}
is obtained in Dawson and Li \cite{DaLi12} and called a
Dawson-Li SDE in Pardoux \cite{Pardoux16}.
SDEs similar to \eqref{1.01} are also studied by
Dawson and Li \cite{DaLi06}
and Fu and Li \cite{FuLi_10}.  

There are a series of results on boundary behaviors to
SDE \eqref{1.01}.
The extinction, extinguishing,
explosion and coming down from infinity are
studied in Li \cite{Li9} with identical power
rate functions $a_0,a_1,a_2$  by Lamperti transform.
Using a martingale approach, the conditions
for extinction, extinguishing,
explosion, coming down from infinity and staying infinite are obtained in
Li et al. \cite{Li et al.10}.
If $a_i(x)=c_i a(x)(i=0,1,2)$ for some constants $c_i$
and function $a(x)$,
then the asymptotic assertions on the speeds of coming down from infinity
and explosion are studied in Foucart et al.
	\cite{Flz19} and in Li and Zhou \cite{Li and Zhou8}, respectively.
The exponential ergodicity is
obtained by Li and Wang \cite{Li and Wang17} by using coupling techniques
and long time behaviour of this SDE with catastrophes
is studied by Marguet and Smadi \cite{Marguet and Smadi14}.
The exponential ergodicity for this SDE with immigration is also studied in \cite{Liliwangzhou}.
The density of the Poisson random measure $N(\dd s,\dd z,\dd u)$
is $\dd s\mu(\dd z)\dd u$, $U=(0,\infty)$ and the condition
$\int_0^\infty (z\wedge z^2)\mu(\dd z)<\infty$
is required in \cite{Li et al.10,Flz19,Li and Zhou8,Marguet and Smadi14}.

Using generalized Chen's criteria techniques, new conditions on non-extinction,
non-explosion, coming down from infinity and staying infinite for the model of
\eqref{1.01} are established in Ma et al. \cite{Ma et al.13}, which gives an affirmative
answer to boundary behavior of critical cases in Li et al. \cite{Li et al.10} and also finds
some non-extinction and non-explosion conditions when
$\int_U(z\wedge z^2)\mu(\dd z)<\infty$ and
$\int_{(0,\infty)\setminus U} (1\vee\ln(1+z))\mu(\dd z)<\infty$.
These techniques are used in Li et al. \cite{Li et al.10}
for the estimation on the first passage probabilities for
SDE \eqref{1.01} with $U=(0,\infty)$
and in Ren et al. \cite{Ren et al.15} for the extinction-extinguishing
behaviors for a two-dimensional general CSBPs.
These criteria generalize similar results for Markov jump processes;
see Chen \cite[Theorems 2.25 and 2.27]{Chen04}.

Inspired by those work,
in this paper we
study the boundary classification 
for a continuous-state nonlinear Neveu's branching  process,
which can be regarded as the solution to the following SDE
that is a modification of \eqref{1.01}:
\begin{align}\label{eq2.1}
X_{t}
=x
+\int_{0}^{t}\int_{0}^{r}\int_{0}^{X_{s-}^{\beta}}z\tilde{M}(\dd s,\dd z,\dd u)
+\int_{0}^{t}\int^{\infty}_{ r}\int_{0}^{X_{s-}^{\theta}}zM(\dd s,\dd z,\dd u),
\end{align}	
where $x,r>0$ and $\beta,\theta\geq0$, $M(\dd s,\dd z,\dd u)$
denotes a Poisson random measure on\ $(0,\infty)^{3}$
with density\ $\dd sz^{-2}\dd z\dd u$, and
$\tilde{M}(\dd s,\dd z,\dd u)=M(\dd s,\dd z,\dd u)-\dd sz^{-2}\dd z\dd u$. 
When $\beta=\theta=1$, the model in \eqref{eq2.1} is considered by Neveu \cite{Neveu22} and the solution
is called continuous-state Neveu's branching process with branching mechanism
\[\psi(u):=\int_0^\infty[\e^{-zu}-1-zu1_{\{z\le r\}}]z^{-2}\dd z=u\ln u +u(\ln r-5\e^{-1}+1).\]
Neveu \cite{Neveu22}  notes a connection between Ruelle's generalized random energy
models  and the continuous-state Neveu's branching processes. 
Bolthausen and Sznitman \cite{BS1998} explains how the results of replica theory of spin glasses can be interpreted in terms
of a coalescent process (now called the Bolthausen-Sznitman coalescent).
Bertoin and Le Gall \cite{bl25} gives a precise and complete form of the relation between the connections between the continuous state
Neveu's branching processes, the Ruelle's generalized random energy
models, and the Bolthausen-Sznitman coalescent.
Further conclusions are obtained in Foucart et al. \cite{FoucartMa}.

Our \textit{aim} is to establish the necessary and sufficient conditions
of the extinction, explosion and coming down from infinity for the continuous-state nonlinear Neveu's branching process determined in \eqref{eq2.1}
and thus find out the sharp influence of different forms of big jumps and small jumps on the boundary behaviors.
To this purpose, we adapt the approach
in Ma et al. \cite{Ma et al.13} for the behavior of non-extinction, non-explosion and
staying infinite and establish new criteria for
the behavior of extinction, explosion and first passage probabilities which are of independent interesting
and generalize the criteria in
\cite[Proposition 2.2]{Ren et al.15} and Chen \cite[Theorem 2.27]{Chen04}
(see Lemma \ref{lemma3.4} in the following).
We also select power functions
and log-logarithm type test functions which are different from that of
Ma et al. \cite{Ma et al.13} (see the explanations in Remark \ref{remark2.1}).

Throughout this paper we always assume
that process\ $(X_{t})_{{t}\geq0}$
is defined on filtered probability space
$(\Omega,\mathscr{F},\mathscr{F}_{t},\mathbf{P})$
which satisfies the usual hypotheses.
We use $\mathbf{P}_{x}$
to denote the law of a process
started at $x$, and denote by $\mathbf{E}_{x}$ the associated expectation.
For $a,b>0$
we define the first passage times
\begin{align}\label{1.2}
\tau_{a}^{-}:=\inf\{t>0:X_{t}\leq a\},~~
\tau_{b}^{+}:=\inf\{t>0:X_{t}\geq b\},
\end{align}
and the first times of
hitting zero $\tau_{0}^{-}$ and reaching infinity
$\tau_{\infty}^{+}$ by
\begin{align*}
\tau_{0}^{-}:=\inf\{t>0:X_{t}=0\} ,~~
\tau_{\infty}^{+}:=\lim_{n\rightarrow\infty}\tau_{n}^{+}
\end{align*}
with the convention $\inf\emptyset=\infty$.
Then $\lim_{a\to0+}\tau_{a}^{-}=\tau_{0}^{-}$.
Let $C^{2}((0,\infty))$ denote the space of twice
continuously differentiable functions on $(0,\infty)$.
We only consider the solution to (\ref{eq2.1})
before $\tau_{0}^{-}\wedge\tau_{\infty}^{+}$,
i.e. both zero and infinity are absorbing states for the
solution.
By the same argument as in Li et al. \cite[Theorem 3.1]{Li et al.10}
and Ren et al. \cite[Lemma A.1]{Ren et al.15},
we can show that SDE \eqref{1.01} has a pathwise unique solution.
Process $(X_{t})_{t\geq 0}$ becomes extinct
iff $\tau_{0}^{-}<\infty$ almost surely; it explodes
iff $\tau_{\infty}^{+}<\infty$ almost surely;
it stays infinite
iff $\lim_{x\rightarrow\infty} \mathbf{P}_{x}\{\tau_{a}^{-}<t\}=0$
for all $t>0$ and all large $a$; it comes down from infinity iff
$\lim_{a\rightarrow\infty}\lim_{ x\rightarrow\infty} \mathbf{P}_x\{\tau_a^-<t\}=1$ for all $t > 0$.
The solution $(X_{t})_{t\geq 0}$ to \eqref{eq2.1} comes down from infinity if it does not stay infinite
by the arguments in \cite[Theorem 1.11]{Li9} and \cite[Lemma 1.2]{Flz20}.
The following assertions are the main result of this paper.

\begin{theorem}\label{theorem2.1}
For the solution $(X_{t})_{t\geq 0}$ to \eqref{eq2.1} we have
\begin{itemize}
\item[(i)]
$(X_{t})_{t\geq0}$ does not extinct iff  $\beta\ge1$;

\item[(ii)]
$(X_{t})_{t\geq0}$ does not explode iff either $\beta>\theta+1$ or $0\le\theta\leq1$;
\item[(iii)]
 $(X_{t})_{t\geq0}$ comes down from infinity iff $\beta>(\theta+1)\vee2$. 
\end{itemize}
\end{theorem}

\begin{remark}\begin{itemize}
\item[(i)]
In the case $\beta=\theta$, the assertions (i) and (ii)  of Theorem \ref{theorem2.1} are given in \cite[Theorems 1.8 and 1.10]{Li9}
but the assertion (iii) of Theorem \ref{theorem2.1} is not presented in \cite{Li9}.
\item[(ii)] From the above example we have the following intuition. The extinction behavior is determined by the 
small jumps term and the process explodes if the big jumps term is much larger than the fluctuation of the 
small jumps term.
The coming down from infinity is due to the fluctuation of the 
small jumps much larger than the big jumps.

\end{itemize}
\end{remark}

The rest of the paper is arranged as follows.
The results on criteria for boundary behaviors are stated in Section 2.
The proof of Theorem \ref{theorem2.1} is deferred to Section 3.

\section{Criteria for boundary behaviors}
\setcounter{equation}{0}

In this section we state the criteria for non-extinction, extinction, non-explosion, explosion and staying infinite,
which generalize Chen's criteria for the uniqueness problem of Markov jump processes.
These Chen's criteria are first proposed in \cite{Chen1986a,Chen1986b} and can also be found in \cite[Theorems 2.25 and 2.27]{Chen04}.

Suppose that $g\in C^{2}((0,\infty))$ satisfies
\begin{align}\label{eq3.2}
\sup_{u\ge v,z\ge1}\Big[|g'(u)|+|g''(u)|+|g(u+z)-g(u)|/z^\delta\Big]<\infty
\end{align}
for any $v>0$ and some $0<\delta<1$. For $u>0$, put
\begin{align}\label{eq3.3}
Lg(u)=u^{\beta}\int_{0}^{r}[g(u+z)-g(u)-zg'(u)]z^{-2}\dd z
 +&u^{\theta}\int^{\infty}_{r}[g(u+z)-g(u)]z^{-2}\dd z.
\end{align}
By\ (\ref{eq2.1}) and It$\rm\hat{o}$'s formula,
we have
 \begin{align}\label{eq3.5}
g(X_{t})=&g(x)+\int_{0}^{t}Lg(X_{s})\dd s+
\int_{0}^{t} \int_{0}^r\int_{0}^{X_{s-}^{\beta}}
[g(X_{s-}+z)-g(X_{s-})]\tilde{M}(\dd s,\dd z,\dd u)\nonumber\\
&+\int_{0}^{t}\int_{r}^\infty\int_{0}^{X_{s-}^{\theta}}
[g(X_{s-}+z)-g(X_{s-})]\tilde{M}(\dd s,\dd z,\dd u).
 \end{align}
For $0<a<b<\infty$
recall the definitions of $\tau_a^-$ and $\tau_b^+$ in \eqref{1.2},
and let $\gamma_{a,b}:=\tau_{a}^{-}\wedge\tau_{b}^{+}$ and
\begin{align*}
M_{t}^{g}:=g(X_{t})-g(x)-\int_{0}^{t}Lg(X_{s})\dd s.
\end{align*}
Then under condition\ (\ref{eq3.2}),
\begin{align}\label{eq3.6}
t\mapsto M_{t\wedge\gamma_{a,b}}^{g}\mbox{ is a martingale}
\end{align}
for all $0<a<b<\infty$.
By Taylor's formula, for any $u,z>0$ we have
\begin{align}\label{eq3.1}
g(u+z)-g(u)-zg'(u)=z^{2}\int_{0}^{1}g''(u+zv)(1-v)\dd v.
\end{align}
The following criteria are
on non-extinction, non-explosion and staying infinite 
for SDE \eqref{eq2.1} and given in \cite[Lemma 3.2]{Ma et al.13}.

\begin{lemma}\label{lemma3.3}
Let $(X_{t})_{t\geq0}$ be the solution to SDE \eqref{eq2.1}.
\begin{itemize}
\item[(i)]
For any fixed $0<b<\infty$, if there exists a function $g\in C^{2}((0,\infty))$ strictly
positive on $(0,b]$ satisfying \eqref{eq3.2} and  $\lim_{u\rightarrow 0}g(u)=\infty$, and
there is a constant $d_b>0$ so that $Lg(u)\leq d_bg(u)$ for all $0<u\le b$,
then $\mathbf{P}_{x}\{\tau_{0}^{-}\geq \tau_{b}^{+}\}=1$ for all $0<x< b$.
\item[(ii)]
For any fixed $0<a<\infty$, if there exists a function $g\in C^{2}((0,\infty))$ strictly
positive on $[a,\infty)$ satisfying \eqref{eq3.2} and $\lim_{u\rightarrow \infty}g(u)=\infty$, and
there is a constant $d_a>0$ so that $Lg(u)\leq d_ag(u)$ for all $u>a$,
then $\mathbf{P}_{x}\{\tau_{\infty}^{+}\geq\tau_{a}^{-}\}=1$ for all $x>a$.
\item[(iii)]
If there exists a function $g\in C^{2}((0,\infty))$ strictly positive on $[u,\infty)$
for all large $u$ satisfying \eqref{eq3.2} and $\lim_{u\rightarrow \infty}g(u)=0$, and
for any large  $a>0$, there is a constant $d_a>0$ such that $Lg(u)\leq d_ag(u)$ for all $u\ge a$,
then $(X_{t})_{t\geq0}$ stays infinite.
\end{itemize}
\end{lemma}

By the arguments in the proof of
\cite[Proposition 2.2]{Ren et al.15} one can get the following
criteria which are on the estimation of 
extinction, explosion and first passage probabilities for SDE \eqref{eq2.1}
as $a=0$, $b=\infty$ and $0<a<b<\infty$ respectively.

\begin{lemma}\label{lemma3.4}
Let $(X_{t})_{t\geq0}$ be the solution to SDE \eqref{eq2.1}.
\begin{itemize}
\item[(i)]
For fixed $0\le a<b<\infty$, if there exists a non-increasing continuous function $g$ on $[0,\infty)$  satisfying $g\in C^{2}((0,\infty))$ and
\eqref{eq3.2}, and there is a
constant $d_{a,b}>0$ so that $Lg(u)\geq d_{a,b}g(u)>0$ for $a<u\le b$,
then $\mathbf{P}_{x}\{\tau_{a}^{-}<\tau_{b}^{+}\}\geq [g(x)-g(b)]/g(a)$ for each $a<x< b$.
\item[(ii)]
For fixed $0<a<b\le \infty$, if there exists a  function $g\in C^{2}((0,\infty))$ satisfying
\eqref{eq3.2} and  $c:=\sup_{a\le u<b} g(u)<\infty$, and there is a
constant $d_{a,b}>0$ so that $Lg(u)\geq d_{a,b}g(u)>0$ for  $a\le u<b$,
then $\mathbf{P}_{x}\{\tau_{b}^{+}<\tau_{a}^{-}\}
\ge [g(x)-g(a)]/c$ for each $a<x<b$.
\end{itemize}
\end{lemma}

\proof
Let  $0<\tilde{a}<\tilde{b}<\infty$. 
Recall \eqref{1.2} and $\gamma_{\tilde{a},\tilde{b}}=\tau_{\tilde{a}}^{-}\wedge\tau_{\tilde{b}}^{+}$.
By (\ref{eq3.5}) and (\ref{eq3.6}), 
for all $g\in C^{2}((0,\infty))$ satisfying \eqref{eq3.2},
\begin{align*}
\mathbf{E}_{x}[g(X_{t\wedge\gamma_{\tilde{a},\tilde{b}}})]=g(x)
+\mathbf{E}_{x}\Big[\int_{0}^{t\wedge\gamma_{\tilde{a},\tilde{b}}}Lg(X_{s})\dd s\Big]
=g(x)+\int_{0}^{t}\mathbf{E}_{x}\Big[Lg(X_{s})1_{\{s\leq\gamma_{\tilde{a},\tilde{b}}\}}\Big]\dd s
\end{align*}
and then by integration by parts, for any $d>0$,
\begin{align}\label{4.17}
&\mathbf{E}_{x}[g(X_{t\wedge\gamma_{\tilde{a},\tilde{b}}})]\e^{-dt} \nonumber\\
=&g(x)-d\int_{0}^{t}\mathbf{E}_{x}[g(X_{s\wedge\gamma_{\tilde{a},\tilde{b}}})\e^{-ds}]\dd s
+\int_{0}^{t}\e^{-ds}\mathbf{E}_{x}[Lg(X_{s})1_{\{s\leq\gamma_{\tilde{a},\tilde{b}}\}}]\dd s.
\end{align}

For part (i), taking $d=d_{a,b}$ in \eqref{4.17}
and using the assumption that $Lg(u)\geq d_{a,b}g(u)>0$
for $a<u\le b$, we have
\begin{align*}
&\mathbf{E}_{x}[g(X_{t\wedge\gamma_{\tilde{a},b}})]\e^{-d_{a,b}t} \nonumber\\
\geq& g(x)-d_{a,b}\int_{0}^{t}\mathbf{E}_{x}[g(X_{s\wedge\gamma_{\tilde{a},b}})]\e^{-d_{a,b}s}\dd s
+d_{a,b}\int_{0}^{t}\e^{-d_{a,b}s}\mathbf{E}_{x}[g(X_{s})1_{\{s\leq\gamma_{\tilde{a},b}\}}]\dd s
\end{align*}
for all $a<\tilde{a}<b$,
which implies that
\begin{align*}
g(x)
\leq\mathbf{E}_{x}[g(X_{t\wedge\gamma_{\tilde{a},b}})]\e^{-d_{a,b}t}
+d_{a,b}\mathbf{E}_{x}\Big[\int_{0}^{t}g(X_{\gamma_{\tilde{a},b}})
\e^{-d_{a,b}s}1_{\{s>\gamma_{\tilde{a},b}\}}\dd s\Big].
\end{align*}
Letting $t\rightarrow\infty$ in the above inequality and using the dominated convergence we obtain
\begin{align*}
g(x)\leq d_{a,b}\mathbf{E}_{x}\Big[g(X_{\gamma_{\tilde{a},b}})\int_{\gamma_{\tilde{a},b}}^{\infty}
\e^{-d_{a,b}s}\dd s\Big]
=\mathbf{E}_{x}[g(X_{\gamma_{\tilde{a},b}})\e^{-d_{a,b}\gamma_{\tilde{a},b}}]
\end{align*}
for all $a<\tilde{a}<b$.
Letting $\tilde{a}\downarrow a$ we have
\begin{align*}
g(x)\leq& \mathbf{E}_{x}\Big[g(X_{\tau_{a}^{-}\wedge\tau_{b}^{+}})
\e^{-d_{a,b}(\tau_{a}^{-}\wedge\tau_{b}^{+})}
(1_{\{\tau_{a}^{-}<\tau_{b}^{+}\}}+1_{\{\tau_{a}^{-}>\tau_{b}^{+}\}}
+1_{\{\tau_{a}^{-}=\tau_{b}^{+}=\infty\}})\Big].
\end{align*}
Since $g(u)\le g(b)$ for all $u\ge b$, then
\begin{align*}
g(x)
\leq g(a)\mathbf{E}_{x}\big[1_{\{\tau_{a}^{-}<\tau_{b}^{+}\}}\big]
+\mathbf{E}_{x}\big[g(X_{\tau_{b}^{+}})\big]
\le g(a)\mathbf{P}_{x}\{\tau_{a}^{-}<\tau_{b}^{+}\}+g(b),
\end{align*}
which gives the assertion (i).

For part (ii), using \eqref{4.17} with $d=d_{a,b}$
and the assumption that $Lg(u)\geq d_{a,b}g(u)>0$
for $a\le u<b$ and by the similar argument one can obtain
\begin{align*}
g(x)
\leq\mathbf{E}_{x}[g(X_{t\wedge\gamma_{a,\tilde{b}}})]\e^{-d_{a,b}t}
+d_{a,b}\mathbf{E}_{x}\Big[\int_{0}^{t}g(X_{\gamma_{a,\tilde{b}}})\e^{-d_{a,b}s}
1_{\{s>\gamma_{a,\tilde{b}}\}}\dd s\Big]
\end{align*}
for all $a<\tilde{b}<b$.
Letting $t\rightarrow\infty$ and then $\tilde{b}\uparrow  b$ we have
\begin{align*}
g(x)\leq& \mathbf{E}_{x}\Big[g(X_{\tau_{a}^{-}\wedge\tau_b^{+}})
\e^{-d_{a,b}(\tau_{a}^{-}\wedge\tau_b^{+})}
(1_{\{\tau_b^{+}<\tau_{a}^{-}\}}
+1_{\{\tau_b^{+}>\tau_{a}^{-}\}})\Big],
\end{align*}
which implies
\begin{align*}
g(x)\leq c\mathbf{P}_{x}\{\tau_b^{+}<\tau_{a}^{-}\}+g(a).
\end{align*}
This shows the assertion (ii).  \qed

\section{Proof of the main result}

\setcounter{equation}{0}

In this section we establish the proof for Theorem \ref{theorem2.1}
by using the generalized Chen's criteria in Lemmas \ref{lemma3.3} and \ref{lemma3.4}. For the assertions of non-extinction and non-explosion
we use Lemma \ref{lemma3.3}(i)-(ii) and select log-logarithm type test functions.
Applying Lemma \ref{lemma3.4} and choosing power functions we establish the behaviors of 
extinction and explosion. We use Lemma \ref{lemma3.3}(iii) choose
a power function for proving the behavior of staying infinite.
We first study the following first passage probabilities.
For the convenience of calculation we always assume that $r=1$ in  \eqref{eq2.1} in this section.

\begin{proposition}\label{proposition2.2}
Let $0<a<x<b<\infty$ be fixed.
For the solution $(X_{t})_{t\geq 0}$ to \eqref{eq2.1} we have
\begin{itemize}
\item[(i)]  $\mathbf{P}_{x}\{\tau_{a}^{-}<\infty\}>0$;
\item[(ii)] $\mathbf{P}_{x}\{\tau_{b}^{+}<\infty\}>0$.
\end{itemize}
\end{proposition}
\proof
For $\lambda,u>0$ define $g_\lambda(u)=\e^{-\lambda u}$.
Then
\begin{align*}
g'_\lambda(u)=-\lambda \e^{-\lambda u},~~g''_\lambda(u)=\lambda^2 \e^{-\lambda u}, \qquad u>0.
\end{align*}
By (\ref{eq3.3}) and \eqref{eq3.1} we have
\begin{align} \label{4.35}
Lg_\lambda(u)=&u^\beta\lambda^2\int^1_0 \dd z\int^1_0\e^{-\lambda (u+zv)}(1-v)\dd v+u^{\theta}\int^\infty_1[\e^{-\lambda (u+z)}-\e^{-\lambda u}]z^{-2}\dd z\nonumber\\
\ge&g_\lambda(u)\Big[2^{-1}u^\beta\lambda^2\int^1_0 \dd z\int^{\frac12}_0\e^{-\lambda z/2}\dd v
-u^{\theta}\int^\infty_1z^{-2}\dd z\Big]\nonumber\\
=&
g_{\lambda}(u)\big[2^{-1}u^{\beta}\lambda (1-\e^{-\frac{\lambda }{2}})-u^{\theta}\big].
\end{align} 
Then there are constants $\lambda_0,d>0$ so that  for all $a\le u\le b$,
\begin{align*}
Lg_{\lambda_0}(u)\geq g_{\lambda_0}(u)\big[2^{-1}a^{\beta}\lambda_0 (1-\e^{-\frac{\lambda_0}{2}})-b^{\theta}\big]
\ge dg_{\lambda_0}(u).
\end{align*}
From Lemma \ref{lemma3.4}(i) it follows that
\begin{align}\label{4.44}
\mathbf{P}_x\{\tau^-_a<\infty\}\geq \mathbf{P}_x\{\tau^-_a<\tau^+_b\}\geq \e^{\lambda_0 a}[\e^{-\lambda_0 x}-\e^{-\lambda_0 b}]>0,
\end{align}
which establishes assertion (i).

Let non-increasing function $g\in C^{2}((0,\infty))$ satisfy $g(u)=\e^{u}$ for $0<u\le b+1$ and $g'(u)=0$ for all $u\ge b+2$.
Then
\begin{align*}
Lg(u)=&u^{\beta}\int_0^1\dd z\int^1_0g''(u+zv)(1-v)\dd v+u^\theta\int_1^\infty [g(u+z)-g(u)]z^{-2}\dd z\nonumber\\
\ge& 2^{-1}u^\beta g(u) \ge 2^{-1}a^\beta g(u),\qquad a\le u\le b.
\end{align*}
By using Lemma \ref{lemma3.4}(ii),
\begin{align}\label{4.43}
\mathbf{P}_x\{\tau^+_b<\infty\}\geq \mathbf{P}_x\{\tau^+_b<\tau^-_a\}\geq (\e^x-\e^a)/\sup_{u>0}g(u)>0,
\end{align}
which completes the proof.\qed

Now we are ready to 

\noindent{\it Proof of Theorem \ref{theorem2.1}(i).}
Suppose that $\beta\ge1$.
For $n\ge1$ let $g_{n}\in C^{2}((0,\infty))$ be
a strictly positive function so that $g_{n}(u)=\ln\ln(6nu^{-1})$ for $0<u\le 2n$ and
\[\sup_{u>2n}[g_{n}(u)+|g_{n}'(u)|+|g_{n}''(u)|]<\infty.\]
Then for each $n\ge1$, there is a constant $c_n>0$ so that
\begin{align}\label{4.0}
g_n(u)\le c_n,\qquad u\ge1.
\end{align}
Moreover,
\begin{align}\label{4.1}
g_n'(u)=-\frac{1}{u\ln(6nu^{-1})},
~~g_n''(u)=\frac{1}{u^{2}\ln(6nu^{-1})}-\frac{1}{u^{2}[\ln(6nu^{-1})]^{2}},
\qquad 0<u\le 2n.
\end{align}
By (\ref{eq3.3}) and a change of variable, we obtain
\begin{align}\label{4.2}
Lg_{n}(u)
=&u^{\beta-1}\int_{0}^{1}[g_{n}(u+uz)-g_{n}(u)-uzg_{n}'(u)]z^{-2}\dd z\nonumber\\
&+u^{\beta-1}\int_{1}^{u^{-1}}[g_{n}(u+uz)-g_{n}(u)-uzg_{n}'(u)]z^{-2}\dd z \nonumber\\
&+u^{\theta-1}\int_{u^{-1}}^{\infty}[g_{n}(u+uz)-g_{n}(u)]z^{-2}\dd u \nonumber\\
=:& u^{\beta-1}[F_{1}(u)+F_{2}(u)]+u^{\theta-1}F_{3}(u),\qquad 0<u\le 1
\end{align}
and
\begin{align}\label{4.3}
Lg_{n}(u)
=&u^{\beta-1}\int_{0}^{u^{-1}}[g_{n}(u+uz)-g_{n}(u)-uzg_{n}'(u)]z^{-2}\dd z\nonumber\\
&+u^{\theta-1}\int_{u^{-1}}^{\infty}[g_{n}(u+uz)-g_{n}(u)]z^{-2}\dd z \nonumber\\
=:&u^{\beta-1}F_{4}(u)+u^{ \theta-1}F_{3}(u),\qquad u>1.
\end{align}
By \eqref{eq3.1} and \eqref{4.1}, for $0<u\le n$ and $0\le z\le1$
we have
\begin{align*}
&g_{n}(u+uz)-g_{n}(u)-uzg_{n}'(u)
=u^{2}z^{2}\int_{0}^{1}g_{n}''(u+uzv)(1-v)\dd v \nonumber\\
&\quad\le u^{2}z^{2}\int_{0}^{1}\frac{1-v}{(u+uzv)^{2}\ln[6n(u+uzv)^{-1}]}\dd v
\leq \frac{z^{2}}{\ln(3nu^{-1})}
\le z^{2}/\ln3,
\end{align*}
which implies
\begin{align}\label{4.4}
F_{1}(u_1)\le (\ln3)^{-1},~F_4(u_2)\le (\ln3)^{-1},
\qquad 0<u_1\le 1,~1\le u_2\le n.
\end{align}
For $0<u\le 1$ and $1\leq z\leq u^{-1}$,
we have $0<u+uz\leq2n$ and then $g_{n}(u+uz)\leq g_{n}(u)$, which implies
\begin{align*}
g_{n}(u+uz)-g_{n}(u)-uzg_{n}'(u)\leq-uzg_{n}'(u)=\frac{z}{\ln(6n u^{-1})},
\qquad 0<u\le 1,~1\leq z\leq u^{-1}
\end{align*}
by \eqref{4.1}. 
It thus follows that
\begin{align}\label{4.5}
F_{2}(u)
\leq\frac{1}{\ln(6nu^{-1})}\int_{1}^{u^{-1}}z^{-1}\dd z
=\frac{\ln u^{-1}}{\ln(6nu^{-1})}\le1,\quad 0<u\le 1.
\end{align}
In view of \eqref{4.0},
\begin{align}\label{4.6}
F_{3}(u)\leq\int_{u^{-1}}^{\infty}g_{n}(u+uz)z^{-2}\dd z
\leq c_n\int_{u^{-1}}^{\infty}z^{-2}\dd z
\le uc_n,\quad 0<u\le n.
\end{align}
Since $\beta\ge 1$, then combining
\eqref{4.4}--\eqref{4.6} with \eqref{4.2}--\eqref{4.3} one obtains
\begin{align}\label{4.7}
Lg_{n}(u)
\le n^{\beta-1}(1+(\ln 3)^{-1})+n^{\theta}c_n,\qquad 0<u\le n.
\end{align}
Observe that $g_n(u)\ge g_n(n)=\ln\ln6$ for all $0<u\le n$.
It then follows from \eqref{4.7} that
 \beqnn
Lg_{n}(u)\le
[(n^{\beta-1}(1+(\ln 3)^{-1})+n^{\theta}c_n)/\ln\ln6]g_n(u),\qquad 0<u\le n.
 \eeqnn
Now using
Lemma~\ref{lemma3.3}(i)
we obtain $\mathbf{P}_{x}\{\tau_{0}^{-}\geq\tau_{n}^{+}\}=1$
for all $0<x< n$.
Letting $n\rightarrow\infty$ we have
$\mathbf{P}_{x}\{\tau_{0}^{-}\geq\tau_\infty^{+}\}=1$
for all $x>0$.
Since the process is defined before the first time of hitting zero or
explosion, $\mathbf{P}_{x}\{\tau_{0}^{-}=\infty$ or $\tau_{\infty}^{+}=\infty\}=1$. Then $\mathbf{P}_{x}\{\tau_{0}^{-}=\infty\}=1$ for all $x>0$.

In the following we assume that $0\le \beta<1$
and $0<\rho<1-\beta$.
Let $0<c<1$ be a small enough constant so that
 \beqlb\label{4.8}
\tilde{c}:=\rho(1-\rho)2^{\rho-3}c^{\beta+\rho-1}-(1-\rho)^{-1}>0.
 \eeqlb
Define the function $g$ on $[0,\infty)$ by $g(u)=c^{\rho}-u^{\rho}$.
Then
\begin{align*}
g'(u)=-\rho u^{\rho-1},~~g''(u)=\rho(1-\rho) u^{\rho-2}, \qquad u>0.
\end{align*}
By (\ref{eq3.3}) and \eqref{eq3.1} we have
\begin{align}\label{4.9}
Lg(u)
=&\rho(1-\rho)u^{\beta}\int_{0}^{1}\dd z\int_{0}^{1}
(u+vz)^{\rho-2}(1-v)\dd v-u^{\theta}\int_{1}^{\infty}[(u+z)^{\rho}-u^{\rho}]z^{-2}\dd z \nonumber\\
=:&G_{1}(u)-G_{2}(u).
\end{align}
By a change of variable we have
\begin{align}\label{4.10}
G_{1}(u)=&\rho(1-\rho)u^{\rho+\beta-1}
\int_{0}^{u^{-1}}\dd z\int_{0}^{1}(1+zv)^{\rho-2}(1-v)\dd v  \nonumber\\
\geq&\rho(1-\rho)u^{\rho+\beta-1}\int_{0}^{1}\dd z
\int_{0}^{1}(1+zv)^{\rho-2}(1-v)\dd v   \nonumber\\
\geq&
\rho(1-\rho)2^{\rho-3}c^{\rho+\beta-1},\qquad 0<u\le c.
\end{align}
Since $(u+z)^{\rho}-u^{\rho}\leq z^{\rho} $ for all $u,z\ge0$,
then
\begin{align}\label{4.11}
G_{2}(u)\leq u^{\theta}\int_{1}^{\infty}z^{\rho-2}\dd z\le(1-\rho)^{-1}.
\end{align}
Combining \eqref{4.8}--\eqref{4.11} we obtain
\begin{align*}
Lg(u)\geq \tilde{c}=\tilde{c}c^{-\rho}c^{\rho}
\ge \tilde{c}c^{-\rho}g(u), \qquad 0<u\leq c.
\end{align*}
By using Lemma~\ref{lemma3.4}(i) one gets that for $0<x<c$,
\begin{align}\label{4.45}
\mathbf{P}_{x}\{\tau_{0}^{-}<\infty\} \ge
\mathbf{P}_{x}\{\tau_{0}^{-}<\tau_{c}^{+}\}
\ge [g(x)-g(c)]/g(0)
=1-(x/c)^{\rho}.
\end{align}
By strong Markov property, \eqref{4.45}  and Proposition \ref{proposition2.2}(i), for all $x\ge c$,
\begin{align*}
\mathbf{P}_{x}\{\tau_{0}^{-}<\infty\}
&\geq\mathbf{P}_{x}\{\tau_{0}^{-}\circ \theta_{\tau_{c/2}^{-}}<\infty,\tau_{c/2}^{-}<\infty\}
\nonumber\\
&=\mathbf{E}_{x}\Big[\mathbf{E}_{x}\big[1_{\{\tau_{0}^{-}\circ \theta_{\tau_{c/2}^{-}}<\infty\}}1_{\{\tau_{c/2}^{-}<\infty\}}|\mathscr{F}_{\tau_{c/2}^-}\big] \Big]\nonumber\\
&=\mathbf{E}_{x}\Big[\mathbf{P}_{c/2} \{\tau_{0}^{-}<\infty\}1_{\{\tau_{c/2}^{-}<\infty\}} \Big]\nonumber\\
&\ge (1-2^{-\rho})\mathbf{P}_{x}\{\tau_{c/2}^{-}<\infty\}>0,
\end{align*}
where $\theta_t$ denotes the usual shift operator. 
This concludes the assertion.
\qed

\noindent{\it Proof of Theorem \ref{theorem2.1}(ii).}
Suppose that  $\beta>\theta+1$ or $\theta\leq1$.
For each $n\geq1$ let $g_{n}\in C^{2}((0,\infty))$ be a strictly positive
function so that $g_{n}(u)=\ln\ln(6nu)$ for $u\ge n^{-1}$.
It then follows that
\begin{align}\label{4.12}
g_{n}'(u)=\frac{1}{u\ln(6nu)},~~
g_{n}''(u)=-\frac{1}{u^{2}\ln(6nu)}-\frac{1}{u^{2}(\ln(6nu))^{2}},
\qquad u\geq n^{-1}.
\end{align}
By (\ref{eq3.3}) and a change of variable, we have
\begin{align}\label{4.13}
Lg_{n}(u)
=&u^{\beta-1}\int_{0}^{\frac{1}{u}}[g_{n}(u+uz)-g_{n}(u)-uzg_{n}'(u)]z^{-2}\dd z
+u^{\theta-1}\int_{\frac{1}{u}}^{\infty}[g_{n}(u+uz)-g_{n}(u)]z^{-2}\dd z\nonumber\\
=:&u^{\beta-1}H_{1}(u)+u^{\theta-1}H_{2}(u).
\end{align}
By \eqref{eq3.1} and \eqref{4.12}, it is elementary to see that
we have
\begin{align*}
g_{n}(u+uz)-g_{n}(u)-uzg_{n}'(u) 
\leq-u^2z^2\int_{0}^{1}\frac{1-v}{(u+uzv)^2\ln(6n(u+uzv))}\dd v
\le -\frac{z^2}{8\ln(12nu)}
\end{align*}
for all $u\ge n^{-1}$ and $0<z\le1$,
which implies
\begin{align}\label{eq3.12}
 H_{1}(u)\leq -\frac{u^{\beta-2}}{8\ln(12nu)},\quad u\ge n^{-1}.
\end{align}
Similarly,
for $u\ge n^{-1}$ and $z\ge0$ we have
\begin{align}\label{eq3.19}
&g_{n}(u+uz)-g_{n}(u)=uz\int_{0}^{1}g_{n}'(u+uzv)\dd v\nonumber\\
&\quad\leq [\ln(6nu)]^{-1}\int_{0}^{1}\frac{z}{1+zv}\dd v=
[\ln(6nu)]^{-1}\ln(1+z).
\end{align}
Observe that
 \beqnn
\lim_{u\to0}\frac{\int_{u}^{\infty}\frac{\ln(1+z)}{z^2}\dd z}{1+\ln u^{-1}}
=\lim_{u\to0}\frac{\frac{\ln(1+u)}{u^2}}{u^{-1}}=1,
 \eeqnn
which implies
 \beqnn
\int_{u^{-1}}^{\infty}\frac{\ln(1+z)}{z^2}\dd z\le
c_0 (1+\ln u),\qquad u\ge1
 \eeqnn
for some constant $c_0>0$.
From (\ref{eq3.19}) it thus follows that
\begin{align}\label{4.14}
H_{2}(u)
\leq[\ln(6nu)]^{-1}\int_{u^{-1}}^{\infty}\frac{\ln(1+z)}{z^2}\dd z
\leq \tilde{c}_0,~~n^{-1}\le u\le1
\end{align}
with $\tilde{c}_0:=(\ln6)^{-1}\int_{1}^{\infty}\frac{\ln(1+z)}{z^{2}}\dd z<\infty$
and
\begin{align}\label{4.15}
H_{2}(u)
\leq[\ln(6nu)]^{-1}\int_{u^{-1}}^{\infty}\frac{\ln(1+z)}{z^2}\dd z
\leq  \frac{c_0(1+\ln u)}{\ln(6n)+\ln u}\le c_0,\qquad u\ge1.
\end{align}
Since $g_n(u)\ge \ln\ln 6$ for all $u\ge n^{-1}$,
then combining \eqref{4.13} with \eqref{eq3.12} and \eqref{4.14}--\eqref{4.15}
one obtains
\begin{align*}
Lg_{n}(u)\le&-u^{\beta-2}[8\ln(12nu)]^{-1}+u^{ \theta-1}(c_0+\tilde{c}_0)\nonumber\\
=&u^{\theta-1} \big[c_0+\tilde{c}_0-u^{\beta-\theta-1}[8\ln(12nu)]^{-1}\big]\nonumber\\
\leq& c_n\le c_n(\ln\ln 6)^{-1}g_n(u),\qquad
u\ge n^{-1}.
\end{align*}
for constants $c_n>0$ as $\beta>\theta+1$.
Similarly, as $0\le\theta\le1$, 
\begin{align*}
Lg_{n}(u)\le u^{\theta-1}(c_0+\tilde{c}_0)
\le n^{1-\theta}(c_0+\tilde{c}_0)(\ln\ln 6)^{-1}g_n(u),\qquad
u\ge n^{-1}.
\end{align*}
Now using
Lemma~\ref{lemma3.3}(ii)
we obtain $\mathbf{P}_{x}\{\tau_{\infty}^{+}\ge\tau_{n^{-1}}^{-}\}=1$
for all $x> n^{-1}$. Letting $n\rightarrow\infty$
we have $\mathbf{P}_{x}\{\tau_{\infty}^{+}\ge\tau_{0}^{-}\}=1$
for all $x>0$. It then follows from the fact of the solution
to SDE (\ref{eq2.1}) defined before $\tau_{0}^{-}\wedge\tau_{\infty}^{+}$ that
$\mathbf{P}_{x}\{\tau_{\infty}^{+}=\infty\}=1$ for all $x>0$.

In the following we assume that $\theta+1\ge\beta$ and $\theta>1$. Let $\rho=\theta-1>0$.
For $u>0$
let $g(u)=1-u^{-\rho}$.
Thus
\begin{align*}
g'(u)=\rho u^{-\rho-1},~g''(u)=-\rho(\rho+1) u^{-\rho-2},
\qquad u>0.
\end{align*}
It then follows from (\ref{eq3.3}) and \eqref{eq3.1} that
\begin{align}\label{4.16}
Lg(u)
=&-\rho(\rho+1)u^{\beta}\int_{0}^{1}\dd z\int_{0}^{1}
(u+vz)^{-\rho-2}(1-v)\dd v+u^{\theta}
\int_{1}^{\infty}[u^{-\rho}-(u+z)^{-\rho}]z^{-2}\dd z \nonumber\\
=:&I_{1}(u)+I_{2}(u).
\end{align}
By a change of variable, for all $u\ge1$,
\begin{align}\label{4.16b}
I_{1}(u)
=&-\rho(1+\rho)u^{\beta-\theta}
\int_{0}^{u^{-1}}\dd z\int_{0}^{1}(1+zv)^{-\rho-2}(1-v)\dd v \nonumber\\
\ge&- \rho(1+\rho)u^{\beta-\theta-1}\ge- \rho(1+\rho)
\end{align}
since $\beta-\theta-1\le0$
and
\begin{align}\label{4.16c}
I_{2}(u)=u^{\theta-\rho-1}\int_{u^{-1}}^{\infty}[1-(1+z)^{-\rho}]z^{-2}\dd z
=\int_{u^{-1}}^{\infty}[1-(1+z)^{-\rho}]z^{-2}\dd z
=:c(u).
\end{align}
Observe that
\begin{align}\label{4.21a}
	\lim_{u\to0}\frac{\int_{u}^{\infty}[1-(1+z)^{-\rho}]z^{-2}\dd z}{\rho\ln u^{-1}}=\lim_{u\to0}
	\frac{[1-(1+u)^{-\rho}]u^{-2}}{\rho u^{-1}}=1,
\end{align}
which implies that $c(u_0)-\rho(1+\rho)>0$ for some large $u_0>1$.
Then combining \eqref{4.16b} and \eqref{4.16c} with \eqref{4.16} we obtain 
\begin{align*}
Lg(u)\ge c(u_0)-\rho(1+\rho)\ge [c(u_0)-\rho(1+\rho)]g(u),\qquad u\ge u_0,
\end{align*}
which gives
\begin{align}\label{4.46}
\mathbf{P}_{x}\{\tau_{\infty}^{+}<\infty\}
\ge\mathbf{P}_{x}\{\tau_{\infty}^{+}<\tau_{u_0}^{-}\}
\ge u_0^{-\rho}-x^{-\rho}
\end{align}
for all $x> u_0$ by Lemma~\ref{lemma3.4}(ii).
By strong Markov property, \eqref{4.46} and Proposition \ref{proposition2.2}(ii), for all $0<x\le u_0$,
\begin{align*}
	\mathbf{P}_{x}\{\tau_{\infty}^{+}<\infty\}
	&\geq\mathbf{P}_{x}\{\tau_{\infty}^{+}\circ \theta_{\tau_{2u_0}^{+}}<\infty,\tau_{2u_0}^{+}<\infty\}
	\nonumber\\
	&=\mathbf{E}_{x}\Big[\mathbf{E}_{x}\big[1_{\{\tau_{\infty}^{+}\circ \theta_{\tau_{2u_0}^{+}}<\infty\}}1_{\{\tau_{2u_0}^{+}<\infty\}}|\mathscr{F}_{\tau_{2u_0}^+}\big] \Big]\nonumber\\
	&=\mathbf{E}_{x}\Big[\mathbf{P}_{2u_0} \{\tau_{\infty}^{+}<\infty\}1_{\{\tau_{2u_0}^{+}<\infty\}} \Big]\nonumber\\
	&\ge u_0^{-\rho} (1-2^{-\rho})\mathbf{P}_{x}\{\tau_{2u_0}^{+}<\infty\}>0,
\end{align*}
which ends the proof.
\qed

	\begin{remark}\label{remark2.1}
At the beginning of the proofs of Theorem \ref{theorem2.1}(i)-(ii),
if the logarithmic type test function is selected as in Ma et al. \cite{Ma et al.13},
the estimations of \eqref{4.5} and \eqref{4.15}  fail to hold.
	\end{remark}

\noindent{\it Proof of Theorem \ref{theorem2.1}(iii).}
If $\beta>\theta+1$ and $\beta\geq2$, then 
$(X_{t})_{t\geq0}$ comes down from infinity
by \cite[Theorem 2.1(iv)]{Ma et al.13} and Theorem \ref{theorem2.1}(ii).
In the following we assume that $\beta\leq(\theta+1)\vee2$. For $\rho>0$ and $u>0$ let $g(u)=u^{-\rho}$.
Then
\begin{align*}
g'(u)=-\rho u^{-\rho-1},~
g''(u)=\rho(\rho+1) u^{-\rho-2},\qquad u>0.
\end{align*}
It then follows from (\ref{eq3.3}) and \eqref{eq3.1} that
\begin{align}\label{4.19}
Lg(u)
=&u^{\beta-1-\rho}
\rho(\rho+1)\int_{0}^{u^{-1}}\dd z\int_0^1(1+zv)^{-\rho-2}(1-v)\dd v\nonumber\\
&-u^{\theta-1-\rho}\int_{u^{-1}}^{\infty}[1-(1+z)^{-\rho}]z^{-2}\dd z , \qquad u>0.
\end{align}
By \eqref{4.21a},
\begin{align}\label{4.21}
\int_{u^{-1}}^{\infty}[1-(1+z)^{-\rho}]z^{-2}\dd z\ge 2^{-1}\rho\ln u,
\qquad u\ge \tilde{c}
\end{align}
for some constant $\tilde{c}>1$.
From (\ref{4.19}) and (\ref{4.21}) it then follows that
\begin{align*}
Lg(u)
\leq u^{\beta-2-\rho}
\rho(\rho+1)
-2^{-1}\rho u^{\theta-1-\rho}\ln u 
=\rho u^{\theta-1-\rho}[(\rho+1)u^{\beta-1-\theta}-2^{-1}\ln u]
\le0,~u\ge c
\end{align*}
for some constant $c\ge \tilde{c}$ as $\beta\le \theta+1$
and 
\begin{align*}
	Lg(u)
	\leq
	\rho(\rho+1) u^{\beta-2-\rho}\le \rho(\rho+1)g(u),\qquad u\ge1
\end{align*}
as $\beta\le2$.
It thus follows that $(X_{t})_{t\geq0}$ stays infinite by Lemma
\ref{lemma3.3}(iii).
\qed

\noindent
{\bf Acknowledgements.}
We are grateful to Pei-Sen Li, Zenghu Li and Xiaowen Zhou for  their careful reading and insightful suggestions on the manuscript.
This work was supported by
NSFC (No. 12061004), NSF of Ningxia (No. 2021AAC02018)
and Major research project for North Minzu University
(No. ZDZX201902).

\end{document}